\documentclass[12pt]{article}
\usepackage{amsmath,amsthm,amsfonts,amssymb}
\usepackage{cite}

\numberwithin{equation}{section}
\theoremstyle{plain}
\newtheorem{theorem}{Theorem}[section]
\newtheorem{corollary}{Corollary}[section]
\newtheorem{proposition}{Proposition}[section]
\newtheorem{lemma}{Lemma}[section]
\theoremstyle{definition}

\theoremstyle{remark}
\newtheorem{remark}{Remark}[section]

\makeatletter
\def\@seccntformat#1{\csname the#1\endcsname.~}
\def\section{\@startsection{section}{1}{\z@}%
              {-3.5ex \@plus -1ex \@minus -.2ex}%
              {2.3ex \@plus.2ex}%
              {\large\bf}}
\def\subsection{\@startsection{subsection}{1}{\z@}%
              {-3.5ex \@plus -1ex \@minus -.2ex}%
              {2.3ex \@plus.2ex}%
              {\normalsize\bf}}

\makeatother

\begin{document}

\title{\large\bf
Transformations and quadratic forms on Wiener spaces
}
\author{\normalsize
Setsuo TANIGUCHI}
\date{}

\maketitle

\begin{abstract}
Two-way relationships between transformations and quadratic
forms on Wiener spaces are investigated with the help of
change of variables formulas on Wiener spaces. 
Further the evaluation of Laplace transforms of quadratic 
forms via Riccati or linear second order ODEs will be
shown. 
\end{abstract}

\renewcommand{\thefootnote}{\arabic{footnote}}
\section{Introduction}

Let $T>0$, $d\in \mathbb{N}$, 
$\mathcal{W}$ be the space of $\mathbb{R}^d$-valued
continuous functions $w$ on $[0,T]$ with $w(0)=0$,
and $\mu$ the Wiener measure on $\mathcal{W}$.
The purpose of this paper is to show two-way relationships
between transformations and quadratic forms on $\mathcal{W}$ 
by use of change of variables formulas on $\mathcal{W}$. 
That is, let $\mathcal{S}_2$ be the space of square
integrable $\mathbb{R}^{d\times d}$-valued
\footnote{
$\mathbb{R}^{d\times d}$ is the space of $d\times d$
real matrices}
functions $\eta=\bigl(\eta_j^i\bigr)_{1\le i,j\le d}$ on
$[0,T]^2$  with $\eta_j^i(t,s)=\eta_i^j(s,t)$ for 
$1\le i,j\le d$ and $(t,s)\in[0,T]^2$.
For $\eta\in \mathcal{S}_2$, define
$G_\eta:\mathcal{W}\to \mathcal{H}$
and $\mathfrak{q}_\eta:\mathcal{W}\to \mathbb{R}$ 
by 
\begin{align*}
    & G_\eta=\biggl(-\sum_{j=1}^d
      \int_0^\bullet \Bigl(\int_0^s
         \eta_j^i(s,u)d\theta^j(u)\Bigr) ds
      \biggr)_{1\le i\le d},
    \\
    & \mathfrak{q}_\eta=\sum_{i,j=1}^d \int_0^T\Bigl(\int_0^t
        \eta_j^i(t,s)d\theta^j(s)\Bigr)d\theta^i(t),
\end{align*}
where 
$\mathcal{H}$ is the Cameron-Martin subspace, 
$\{\theta(t)=(\theta^1(t),\dots,\theta^d(t))
        \}_{t\in[0,T]}$ 
is the coordinate process of $\mathcal{W}$, that is,
$\theta(t)(w)=w(t)$ for $t\in[0,T]$ and  
$w\in \mathcal{W}$, and
$d\theta^i(t)$ is the It\^o integral
with respect to $\{\theta^i(t)\}_{t\in[0,T]}$.
Setting 
\[
    \mathcal{T}=\{G_\eta;\eta\in \mathcal{S}_2\}
    \quad\text{and}\quad
    \mathcal{Q}=\{\mathfrak{q}_\rho;
        \rho\in \mathcal{S}_2\},
\]
we shall show the two-way relationship between $\mathcal{T}$ 
and $\mathcal{Q}$ obtained through the identity
\begin{equation}\label{eq.cv0}
    \int_{\mathcal{W}} f(\iota+G_\eta)
          e^{\mathfrak{q}_\rho} \,d\mu
    =e^{\|\eta\|_2^2/4}
    \int_{\mathcal{W}} f d\mu
    \quad f\in C_b(\mathcal{W}),
\footnote{
$C_b(\mathcal{W})$ is the
space of bounded and continuous functions on $\mathcal{W}$
with values in $\mathbb{R}$.
}
\end{equation}
where $\iota:\mathcal{W}\to \mathcal{W}$ is the identity
map, and
\[
    \|\eta\|_2
    =\Bigl(\int_0^T \int_0^T |\eta(t,s)|^2 dsdt
      \Bigr)^{1/2},
\]
$|M|$ being the Euclidean norm of 
$M\in \mathbb{R}^{d\times d}$.
We shall first present the way from
$\mathcal{T}$ to $\mathcal{Q}$ (constructing
$\mathfrak{q}_\rho$ from given $G_\eta$), and next the
converse way from $\mathcal{Q}$ to $\mathcal{T}$ (showing
the existence of $G_\eta$ producing given
$\mathfrak{q}_\rho$).
See Theorems~\ref{thm.t.q} and \ref{thm.rho.eta}.

It is known that every element of $\mathcal{C}_2$, the
Wiener chaos of order two, is of the form
$\mathfrak{q}_\rho$ as described above. 
For example, see \cite{n}. 
A lot of studies of Wiener integrals 
$\int_{\mathcal{W}} e^{\mathfrak{q}} f d\mu$, 
$\mathfrak{q}\in \mathcal{C}_2$,
are made from various stochastically analytic points of view 
(\cite{cm,ikm2,im1,ikt,it,kac,levy,t,yor} and references
therein). 
The identity \eqref{eq.cv0} gives a method of evaluating
the Wiener integrals 
$\int_{\mathcal{W}} e^{\mathfrak{q}} f d\mu$, 
$\mathfrak{q}\in \mathcal{C}_2$.
Based on the fact that each $\mathfrak{q}\in \mathcal{C}_2$
is specified by a symmetric Hilbert-Schmidt operator $B$
from $\mathcal{H}$ to itself (cf. Remark~\ref{rem.chaos}
below), Malliavin and the author \cite{mall-t} achieved the
identity slightly different from 
\eqref{eq.cv0}:
\begin{equation}\label{eq.mt}
    \int_{\mathcal{W}} e^{\mathfrak{q}} f d\mu
    =\{{\det}_2(I-B)\}^{-1/2} \int_{\mathcal{W}}
       [f\circ(\iota+J_B)]d\mu,
\end{equation}
where $\det_2$ is the regularized determinant,
$I:\mathcal{H}\to \mathcal{H}$ is the identity map, and
$J_B$ is an $\mathcal{H}$-valued random variable obtained
from $B$.
To show this identity, there are two key ingredients: one is
also the change of variables formula on $\mathcal{W}$ via 
Malliavin calculus (\cite{mt}, or Lemma~\ref{lem.cvf}
below); the other is the development of Wiener process
$\{\theta(t)\}_{t\in[0,T]}$ by the ONB of $\mathcal{H}$
consisting of eigenfunctions of $B$.
Applying change of variables formulas to the integrals
$\int_{\mathcal{W}} e^{\mathfrak{q}} f d\mu$ goes back to
Cameron-Martin \cite{cm,cm1} and applying developments of
Wiener processes does to Kac \cite{kac} and 
L\'evy \cite{levy}.

Given $\mathfrak{q}=\mathfrak{q}_\rho$, our above result via
the identity \eqref{eq.cv0} guarantees only the existence of
the corresponding $G_\eta$. 
Further, $\det_2(I-B)$ in \eqref{eq.mt} is rather abstract.
Thus it is natural to ask if there is a more computable
evaluation of $\int_{\mathcal{W}} e^{\mathfrak{q}} f d\mu$.
The second aim of this paper is to meet such demands in a
concrete case. 
That is, letting $\mathcal{T}_0$ be the totality of all
$G_\eta$ with $\eta$ of the form $\eta(t,s)=\chi(t)$, $t>s$,
for some $\chi\in C([0,T];\mathbb{R}^{d\times d})$
\footnote{
$C([0,T];\mathbb{R}^{d\times d})$ is the space of 
continuous $\mathbb{R}^{d\times d}$-valued functions on
$[0,T]$.
}
and $\mathcal{Q}_0$ be that of all $\mathfrak{q}_\rho$
with $\rho$ of the form $\rho(t,s)=\sigma(t)$, $t>s$, for 
some $\sigma\in C([0,T];\mathbb{R}^{d\times d})$, we shall
show the two-way relationship between $\mathcal{T}_0$ and
$\mathcal{Q}_0$. 
Precisely speaking, as for the way from $\mathcal{T}_0$ to
$\mathcal{Q}_0$, for given $G_\eta\in \mathcal{T}_0$, 
we shall show that 
$(\iota+G_\eta)^{-1}:\mathcal{W}\to \mathcal{W}$ exists and 
there is a $\mathfrak{q}_\rho\in \mathcal{Q}_0$ satisfying
the following identity instead of  
\eqref{eq.cv0}:
\begin{equation}\label{eq.cv0+}
    \int_{\mathcal{W}} f e^{\mathfrak{q}_\rho} \,d\mu
    =e^{\int_0^T [\text{\rm tr}\,(\chi(t)-\sigma(t))]dt/2}
    \int_{\mathcal{W}} [f\circ(\iota+G_\eta)^{-1}] d\mu
    \quad f\in C_b(\mathcal{W}).
\end{equation}
See Theorem~\ref{thm.TtoQ}.
To show the way from $\mathcal{Q}_0$ to $\mathcal{T}_0$,
ODEs play a key role. 
In fact, given $\sigma\in C([0,T];\mathbb{R}^{d\times d})$,
if there exists $S\in C^1([0,T];\mathbb{R}^{d\times d})$
\footnote{
$C^k([0,T];\mathbb{R}^{d\times d})$ is the space of
  $k$-times continuously differentiable
  $\mathbb{R}^{d\times d}$-valued functions on $[0,T]$.
}
obeying the matrix Riccati ODE 
\[
    S^\prime=-S^2-\sigma^\dagger S-S\sigma
      -\sigma^\dagger\sigma,
    \quad S(T)=0,
\]
where the symbol ${}^\prime$ means taking
differentiation in $t\in[0,T]$, 
and $S^\dagger(t)=S(t)^\dagger$
\footnote{
$M^\dagger$ is the transpose of $M\in \mathbb{R}^{d\times d}$.
}
for $t\in[0,T]$, 
then $G_\eta\in \mathcal{T}_0$ with $\chi=S+\sigma$ is the
desired transformation corresponding to
$\mathfrak{q}_\rho\in \mathcal{Q}_0$ via
\eqref{eq.cv0+}. 
Further, if $\sigma$ is continuously differentiable
and the solution  
$A\in C^2([0,T];\mathbb{R}^{d\times d})$ to the second order
ODE on $\mathbb{R}^{d\times d}$
\[
    A^{\prime\prime}-2\sigma_A A^\prime-\sigma^\prime A=0,
    \quad A(T)=I_d,~A^\prime(T)=\sigma(T),
\]
where $\sigma_A=\frac12(\sigma-\sigma^\dagger)$ and 
$I_d$ is the $d\times d$ identity matrix, is 
non-singular, that is, $\det A(t)\ne0$ for any
$t\in[0,T]$, 
then $G_\eta$ with $\chi=A^\prime A^{-1}$ is the
transformation corresponding to $\mathfrak{q}_\rho$.
In both cases, the factor 
$e^{{\int_0^T [\text{\rm tr}\,(\chi(t)-\sigma(t))]dt/2}}$ 
in \eqref{eq.cv0+} is given more explicitly in terms of $S$
or $A$.  
See Theorem~\ref{thm.QtoT}.

The two-way relationship between $\mathcal{T}$ and
$\mathcal{Q}$ will be seen in Section~\ref{sec.cv}.
The special cases when $\eta(t,s)=\chi(t)$, $t>s$, i.e., 
the two-way relationship between $\mathcal{T}_0$ and
$\mathcal{Q}_0$ will be
investigated in Section~\ref{sec.TtoQ}. 
In the section, an explicit expression of
$(\iota+G_\eta)^{-1}$ will be given and applied to compute
the conditional expectation of $e^{\mathfrak{q}_\rho}$.
At the end of the same section, two applications of the way
from $\mathcal{Q}_0$ to $\mathcal{T}_0$ will be presented.

\section{General transformations}
\label{sec.cv}

In this section, we shall show the two-way relationship
between the classes $\mathcal{T}$ and $\mathcal{Q}$.

Recall that the Cameron-Marin subspace $\mathcal{H}$
consists of absolutely continuous $h\in \mathcal{W}$ with
the square integrable derivative $h^\prime$, and
it is a real separable Hilbert space equipped with the inner 
product 
\[
    \langle h,g\rangle_{\mathcal{H}}
    =\int_0^T \langle h^\prime(t),g^\prime(t)\rangle dt,
    \quad h,g\in \mathcal{H},
\]
where $\langle\cdot,\cdot\rangle$ is the inner product of
$\mathbb{R}^d$.
Further, remember that the space $\mathcal{S}_2$ was defined
as  
\[
     \mathcal{S}_2
      =\biggl\{\eta:[0,T]^2\to \mathbb{R}^{d\times d};
       \|\eta\|_2<\infty\text{ and }
        \eta(t,s)^\dagger=\eta(s,t),
        ~(t,s)\in [0,T]^2
       \biggr\}.
\]
In what follows, for the sake of simplicity of notation, 
we use the matrix notation;
each element of $\mathbb{R}^d$ is thought of as a column
vector and $\mathbb{R}^{d\times d}$ acts on $\mathbb{R}^d$
from left. 
In particular, the transformation 
$G_\eta:\mathcal{W}\to \mathcal{H}$ and the Wiener
functional $\mathfrak{q}_\eta:\mathcal{W}\to \mathbb{R}$ for  
$\eta=(\eta_j^i)_{1\le i,j\le d}\in \mathcal{S}_2$, which
were given in the previous section, are represented as
\begin{align}
    & \langle G_\eta, h\rangle_{\mathcal{H}}
       =-\int_0^T \Bigl\langle\int_0^t 
           \eta(t,s)d\theta(s),h^\prime(t)\Bigr\rangle
         dt, 
      \quad  h\in \mathcal{H},
 \label{eq.F_eta}
    \\
    & \mathfrak{q}_\eta=\int_0^T \Bigl\langle\int_0^t 
           \eta(t,s)d\theta(s), d\theta(t)
           \Bigr\rangle.
 \label{eq.q_eta}
\end{align}
The first aim of this section is to show the way from
$\mathcal{T}$ to $\mathcal{Q}$.

\begin{theorem}\label{thm.t.q}
Let $\eta\in \mathcal{S}_2$.
Suppose that $\|\eta\|_2<1$.
Take $\rho\in\mathcal{S}_2$ such that
\begin{equation}\label{eq.rho}
    \rho(t,s)=\eta(t,s)-\int_t^T 
          \eta(t,u)\eta(u,s) du
    \quad\text{for }0\le s<t\le T.
\end{equation}
Then \eqref{eq.cv0} holds:
\[
    \int_{\mathcal{W}} f(\iota+G_\eta) e^{\mathfrak{q}_\rho} d\mu
    =e^{\|\eta\|_2^2/4} \int_{\mathcal{W}} fd\mu,
    \quad f\in C_b(\mathcal{W}).
\]
\end{theorem}

Notice that $\|\rho\|_2<\infty$.
In fact, by the Schwarz inequality, it holds that
\[
    \|\rho-\eta\|_2^2
    =2\int_0^T \Bigl(\int_0^t \Bigl|\int_t^T 
     \eta(t,u)\eta(u,s) du \Bigr|^2 ds\Bigr) dt
    \le \frac12 \|\eta\|_2^4.
\]
The proof of Theorem~\ref{thm.t.q}  will be broken into
several steps, each being a lemma.

For a real separable Hilbert space $E$, let
$\mathbb{D}^\infty(E)$ be the space of infinitely
$\mathcal{H}$-differen\-tiable Wiener functionals in the
sense of Malliavin calculus, whose $\mathcal{H}$-derivatives
of all orders are $p$th integrable with respect to $\mu$ for
every $p\in(1,\infty)$. 
The $\mathcal{H}$-derivative and its adjoint are written by
$D$ and $D^*$, respectively.
Both 
$D:\mathbb{D}^\infty(E)\to
   \mathbb{D}^\infty(\mathcal{H}\otimes E)$ and
$D^*:\mathbb{D}^\infty(\mathcal{H}\otimes E)\to
   \mathbb{D}^\infty(E)$ are continuous,
where $\mathcal{H}\otimes E$ is the Hilbert space of
Hilbert-Schmidt operators from $\mathcal{H}$ to $E$.
For details, see \cite{mt}.

Regarding a symmetric 
$B\in \mathcal{H}^{\otimes 2}
 =\mathcal{H}\otimes \mathcal{H}$ 
as a constant function belonging to 
$\mathbb{D}^\infty(\mathcal{H}^{\otimes 2})$, define the
Wiener functional $Q_B\in \mathbb{D}^\infty(\mathbb{R})$ by 
\[
    Q_B=(D^*)^2 B,
\]
and call it the {\it quadratic form} associated with $B$.
The reason why it is called so can be seen in the
following assertion.

\begin{lemma}\label{lem.quad}
If $G\in \mathbb{D}^\infty(\mathbb{R})$ satisfies that
$D^3G=0$, then $D^2G$ is a constant, 
say $B\in \mathcal{H}^{\otimes 2}$, and it holds that
\[
    G=c+D^*h+\frac12 Q_B,
    \quad\text{with }
    c=\int_{\mathcal{W}} G d\mu \text{ and }
    h=\int_{\mathcal{W}} DG d\mu.
\]
Conversely, for any symmetric 
$B\in \mathcal{H}^{\otimes 2}$, it holds that
$D^3Q_B=0$, $\int_{\mathcal{W}}Q_Bd\mu=0$, and
$\int_{\mathcal{W}} DQ_B d\mu=0$.
\end{lemma}

\begin{proof}
See \cite[Propositions~5.2.9 and 5.7.4]{mt}.
\end{proof}

In what follows, we fix $\eta\in\mathcal{S}_2$.
Define $B_\eta:\mathcal{H}\to \mathcal{H}$ by
\begin{equation}\label{eq.B_eta}
    (B_\eta h)^\prime(t)
    =\int_0^T \eta(t,s)h^\prime(s) ds,
    \quad t\in[0,T],~h\in \mathcal{H}.
\end{equation}

\begin{lemma}\label{lem.B_eta}
$B_\eta$ is a symmetric Hilbert-Schmidt operator, 
and satisfies that
\begin{equation}\label{l.B_eta.1}
    \mathfrak{q}_\eta=\frac12 Q_{B_\eta}.
\end{equation}
Further, 
$e^{\lambda\mathfrak{q}_\eta}\in \bigcup_{p\in(1,\infty)}L^p(\mu)$ 
\footnote{
$L^p(\mu)$ is the
space of $p$th integrable $\mathbb{R}$-valued Wiener
functionals with respect to $\mu$.}
for $\lambda\in \mathbb{R}$ with
$|\lambda|\|\eta\|_2<1$.
\end{lemma}

\begin{remark}\label{rem.chaos}
Every $G\in \mathcal{C}_2$ admits $\eta\in\mathcal{S}_2$
such that $G=\mathfrak{q}_\eta$.
For example, see \cite{n}.
Moreover, by this lemma, defining the symmetric
$B\in\mathcal{H}^{\otimes 2}$ by $B=D^2G$, we have that
$G=Q_B/2$. 
\end{remark}

\begin{proof}
It is an easy exercise of Malliavin calculus (cf.\cite{mt})
to see that 
\begin{align*}
    \langle D\mathfrak{q}_\eta,h\rangle_{\mathcal{H}}
      & =\int_0^T \Bigl\langle\int_0^t 
        \eta(t,s)h^\prime(s)ds, d\theta(t)
        \Bigr\rangle
    \\
    & \qquad\qquad
      +\int_0^T \Bigl\langle\int_0^t 
        \eta(t,s)d\theta(s), h^\prime(t)
        \Bigr\rangle dt,
    \\
    \langle (D^2\mathfrak{q}_\eta)[g],h\rangle_{\mathcal{H}}
    & =\int_0^T \Bigl\langle\int_0^t 
        \eta(t,s)h^\prime(s) ds, g^\prime(t)
        \Bigr\rangle dt
    \\
    & \qquad\qquad
      +\int_0^T \Bigl\langle\int_0^t 
        \eta(t,s)g^\prime(s) ds, h^\prime(t)
        \Bigr\rangle dt,
\end{align*}
where $(D^2\mathfrak{q}_\eta)[g]$ is the Wiener functional
whose value at $w\in \mathcal{W}$ is the value of the
Hilbert-Schmidt operator $(D^2\mathfrak{q}_\eta)(w)$ at
$g\in \mathcal{H}$. 
Changing the order of integration and using the relation
that $\eta(t,s)^\dagger=\eta(s,t)$, we see that
\[
    \int_0^T \Bigl\langle\int_0^t 
        \eta(t,s)h^\prime(s)ds, g^\prime(t)
        \Bigr\rangle dt
    =\int_0^T \Bigl\langle\int_t^T 
        \eta(t,s) g^\prime(s)ds, h^\prime(t)
        \Bigr\rangle dt.
\]
Thus $D^2\mathfrak{q}_\eta=B_\eta$, which also implies that 
$B_\eta$ is a symmetric Hilbert-Schmidt operator.
By the above identities, we have that
$D^3\mathfrak{q}_\eta=0$ and  
$\int_{\mathcal{W}}D\mathfrak{q}_\eta d\mu=0$.
Moreover, it is easily seen that 
$\int_{\mathcal{W}}\mathfrak{q}_\eta d\mu=0$.
Due to Lemma~\ref{lem.quad}, \eqref{l.B_eta.1} holds.

It was seen in \cite[Example~5.4.3]{mt} that, for symmetric
$B\in \mathcal{H}^{\otimes 2}$, 
$e^{|\lambda||Q_B|}\in L^1(\mu)$ for $\lambda\in \mathbb{R}$
with $|\lambda|\|B\|_{\text{op}}<\frac12$, where
$\|B\|_{\text{op}}$ is the operator norm of $B$.
By using the Schwarz inequality, it is easily seen that
$\|B_\eta\|_{\text{op}}\le\|\eta\|_2$.
Hence the proof of the second assertion completes.
\end{proof}

The change of variables formula on $\mathcal{W}$, which we
shall use, is stated as follows.

\begin{lemma}\label{lem.cvf}
Let $G\in \mathbb{D}^\infty(\mathcal{H})$.
Suppose that there exists $r\in(\frac12,\infty)$ such
that  
\[
    e^{-D^*G+r\|DG\|_{\mathcal{H}^\otimes 2}^2}
    \in \bigcup_{p\in(1,\infty)} L^p(\mu),
\]
where 
$\|\cdot\|_{\mathcal{H}^{\otimes2}}$ is the Hilbert norm of
$\mathcal{H}^{\otimes 2}$.
Then it holds that
\[
    \int_{\mathcal{W}} f(\iota+G) {\det}_2(I+DG)
      e^{-D^*G-\frac12\|G\|_{\mathcal{H}}^2} d\mu
      =\int_{\mathcal{W}} f d\mu, 
    \quad
    f\in C_b(\mathcal{W}).
\]
\end{lemma}

\begin{proof}
See \cite[Theorem~5.6.1]{mt}.
\end{proof}

We apply this lemma to $G_\eta$ defined in 
\eqref{eq.F_eta}.

\begin{lemma}\label{lem.cvf2}
Suppose that $\|\eta\|_2<1$.
Then it holds that
\begin{equation}\label{l.cvf2.1}
    \int_{\mathcal{W}} f(\iota+G_\eta) 
        e^{\mathfrak{q}_\eta-\mathfrak{h}_\eta} d\mu
    =\int_{\mathcal{W}} f d\mu,
    \quad f\in C_b(\mathcal{W}),
\end{equation}
where
\begin{equation}\label{l.cvf2.2}
    \mathfrak{h}_\eta=\frac12\int_0^T 
       \Bigl| \int_0^t  \eta(t,s)d\theta(s) \Bigr|^2
       dt.
\end{equation}
\end{lemma}

\begin{proof}
By \cite[Theorem~5.3.3]{mt}, we see that
\begin{equation}\label{l.cvf2.21}
    D^*G_\eta=-\mathfrak{q}_\eta.
\end{equation}

Taking the $\mathcal{H}$-derivatives of both sides of
\eqref{eq.F_eta}, we have that 
$\mathcal{H}^{\otimes 2}$-valued Wiener functional $DG_\eta$
satisfies that 
\[
    \langle (DG_\eta)[g],h\rangle_{\mathcal{H}}
    =-\int_0^T \Bigl\langle\int_0^t 
      \eta(t,s)g^\prime(s)ds,h^\prime(t) 
        \Bigr\rangle dt,
    \quad g,h\in \mathcal{H}.
\]
This means that
\begin{equation}\label{l.cvf2.22}
    ((DG_\eta)[g])^\prime(t)
    =-\int_0^t \eta(t,s)g^\prime(s) ds, 
    \quad g\in \mathcal{H}.
\end{equation}
Thus $DG_\eta$ is a Volterra operator, and hence
\begin{equation}\label{l.cvf2.23}
    {\det}_2(I+DG_\eta)=1.
\end{equation}

Letting $\{e_n\}_{n=1}^\infty$ be an ONB of $\mathcal{H}$, 
by \eqref{l.cvf2.22}, we have that
\[
    \|DG_\eta\|_{\mathcal{H}^{\otimes 2}}^2
    =\sum_{n=1}^\infty \int_0^T 
        \Bigl|\int_0^t \eta(t,s)e_n^\prime(s) ds\Bigr|^2
        dt
     =\frac12 \|\eta\|_2^2.
\]
In conjunction with \eqref{l.cvf2.21} and
Lemma~\ref{lem.B_eta}, this yields that
\[
    e^{-D^*G_\eta+r\|DG_\eta\|_{\mathcal{H}^{\otimes 2}}^2}
    \in \bigcup_{p\in(1,\infty)} L^p(\mu)
    \quad\text{for any }r\in[0,\infty).
\]
Since $\frac12\|G_\eta\|_{\mathcal{H}}^2=\mathfrak{h}_\eta$,
applying Lemma~\ref{lem.cvf} to $G=G_\eta$ with use of
\eqref{l.cvf2.21} and \eqref{l.cvf2.23}, we arrive at
\eqref{l.cvf2.1}. 
\end{proof}

\begin{lemma}\label{lem.C_eta}
Define $C_\eta:\mathcal{H}\to \mathcal{H}$ by
\[
    (C_\eta g)^\prime(t)=\int_t^T \Bigl(\int_0^s 
     \eta(t,s) \eta(s,u) g^\prime(u) du\Bigr)ds,
    \quad t\in[0,T],~g\in \mathcal{H}.
\]
Then $C_\eta$ is a symmetric Hilbert-Schmidt operator, 
and satisfies that
\begin{equation}\label{l.C_eta.1}
    \mathfrak{h}_\eta=\frac12 Q_{C_\eta}+\frac14 \|\eta\|_2^2.
\end{equation}
\end{lemma}

\begin{proof}
It is easily seen that
\begin{equation}\label{l.C_eta.21}
    \int_{\mathcal{W}} \mathfrak{h}_\eta d\mu
    =\frac14\|\eta\|_2^2.
\end{equation}

Observe that
\begin{equation}\label{l.C_eta.22}
    \langle D\mathfrak{h}_\eta,h\rangle_{\mathcal{H}}
    =\int_0^T \Bigl\langle
        \int_0^t \eta(t,s)h^\prime(s) ds,
        \int_0^t \eta(t,s)d\theta(s)
     \Bigl\rangle dt,
   \quad h\in \mathcal{H}. 
\end{equation}
This implies that
\begin{equation}\label{l.C_eta.23}
    \int_{\mathcal{W}} D\mathfrak{h}_\eta d\mu=0.
\end{equation}

Taking the $\mathcal{H}$-derivatives of both sides of
\eqref{l.C_eta.22}, we obtain that
\begin{align*}
    \langle (D^2\mathfrak{h}_\eta)[g],h\rangle_{\mathcal{H}}
    & =\int_0^T \Bigl\langle 
      \int_0^t \eta(t,s)h^\prime(s) ds,
      \int_0^t \eta(t,u)g^\prime(u) du
      \Bigr\rangle dt 
    \\
    & =\int_0^T \Bigl(\int_s^T \Bigl(\int_0^t 
        \bigl\langle 
        \eta(s,t) \eta(t,u) g^\prime(u),
        h^\prime(s) \bigr\rangle du\Bigr)dt\Bigr)ds
    \\
    & =\int_0^T \langle (C_\eta g)^\prime(s),h^\prime(s)
        \rangle ds
\end{align*}
for $g,h\in \mathcal{H}$.
Thus it holds that $D^2\mathfrak{h}_\eta=C_\eta$.
Hence $C_\eta$ is a symmetric Hilbert-Schmidt
operator.
Moreover, by Lemma~\ref{lem.quad} with this,
\eqref{l.C_eta.21}, and \eqref{l.C_eta.23}, we obtain
\eqref{l.C_eta.1}. 
\end{proof}

\begin{lemma}\label{lem.B_rho}
Take $\rho\in\mathcal{S}_2$ satisfying \eqref{eq.rho}.
Then $B_\eta-C_\eta=B_\rho$ and it holds that
\[
    \mathfrak{q}_\eta-\mathfrak{h}_\eta
    =\mathfrak{q}_\rho-\frac14\|\eta\|_2^2.
\]
\end{lemma}

\begin{proof}
Observe the representation
\[
    (C_\eta g)^\prime(t)
    =\int_0^T \int_0^T
     \boldsymbol{1}_{(t,T]}(s)\boldsymbol{1}_{[0,s)}(u)
     \eta(t,s)\eta(s,u)g^\prime(u) du ds,
\]
where $\boldsymbol{1}_A$ is the indicator function of $A$.
Since
\[
     \boldsymbol{1}_{(t,T]}(s)\boldsymbol{1}_{[0,s)}(u)
     =\boldsymbol{1}_{(t,T]}(s)\boldsymbol{1}_{[0,t]}(u)
      +\boldsymbol{1}_{(t,T]}(u)\boldsymbol{1}_{(u,T]}(s),
\]
changing the order of integration, we obtain that
\[
     (C_\eta g)^\prime(t)
     =\int_0^T \Bigl(\int_{t\vee u}^T 
          \eta(t,s)\eta(s,u) ds\Bigr)g^\prime(u) du,
    \quad t\in[0,T],g\in \mathcal{H},
\]
where $t\vee u=\max\{t,u\}$.
Being in $\mathcal{S}_2$, $\rho$ satisfies that
\[
    \rho(t,s)=\eta(t,s)-\int_{s\vee t}^T 
       \eta(t,u)\eta(u,s) du \quad\text{for }t\ne s.
\]
By this, the above expression of $C_\eta g$, and 
the definition \eqref{eq.B_eta} of $B_\eta$, we have that
$B_\eta-C_\eta=B_\rho$.
Then, by Lemmas~\ref{lem.B_eta} and \ref{lem.C_eta}, we have
that 
\[
    \mathfrak{q}_\rho=\frac12 Q_{B_\rho}
    =\frac12 Q_{B_\eta}-\frac12 Q_{C_\eta}
    =\mathfrak{q}_\eta-\mathfrak{h}_\eta+\frac14\|\eta\|_2^2,
\]
which implies the desired identity.
\end{proof}

\begin{lemma}
The assertion of Theorem~\ref{thm.t.q} holds.
\end{lemma}

\begin{proof}
This follows from Lemmas~\ref{lem.cvf2} and \ref{lem.B_rho}.
\end{proof}

We next see the way from $\cal{Q}$ to $\mathcal{T}$.

\begin{theorem}\label{thm.rho.eta}
There is $\varepsilon>0$ such that 
each $\rho\in \mathcal{S}_2$ with $\|\rho\|_2<\varepsilon$
admits $\eta\in \mathcal{S}_2$ such that
$\|\eta\|_2<1$ and the identities \eqref{eq.rho} and 
\eqref{eq.cv0} hold.
\end{theorem}

\begin{proof}
Let $\rho\in \mathcal{S}_2$.
Suppose that $\|\rho\|_2<\frac13$.
Define $\rho^{*n}$ by 
\[
    \rho^{*1}=\rho,\quad
    \rho^{*n}(t,s)=\int_0^T \rho(t,u)\rho^{*(n-1)}(u,s)du,
    \quad n\ge2.
\]
Since $\rho^{*n}\in \mathcal{S}_2$ and 
$\|\rho^{*n}\|_2\le \|\rho\|_2^n$, $n\in \mathbb{N}$, 
the function $\varphi$ defined by
\[
    \varphi=\sum_{n=1}^\infty \rho^{*n}
\]
is in $\mathcal{S}_2$ and $\|\varphi\|_2<\frac12$.

For $\psi\in \mathcal{S}_2$, define the bounded linear
operator $K_\psi$ from $L^2([0,T];\mathbb{R}^d)$
\footnote{
$L^2([0,T];\mathbb{R}^d)$ is the space of square integrable
$\mathbb{R}^d$-valued functions on $[0,T]$ with respect to
the Lebesgue measure.
}
to itself by
\[
    (K_\psi f)(t)=\int_0^T \psi(t,s)f(s)ds,
    \quad f\in L^2([0,T];\mathbb{R}^d).
\]
Then $(I+K_\varphi)^{-1}$, where $I$ is the
identity map of $L^2([0,T];\mathbb{R}^d)$, is a bounded
linear operator and it holds that 
\[
    (I+K_\varphi)^{-1}-I=K_{-\rho}.
\]
Hence $-\rho$ is the resolvent kernel of $-K_\varphi$.
Thanks to the special factorization of $-K_\varphi$ due to
Gohberg and Krein \cite{gk}, there exists
$v\in \mathcal{S}_2$ such that 
\[
    -\rho(t,s)=v(t,s)+\int_t^T v(t,u)v(u,s)du,
    \quad s<t.
\]
This fact is obtained by combining the observations in
\cite{gk} 
(the proposition $1^\circ$ before Theorem~6.2, the
identity (8.5), the remark after (2.5), and Theorem~3.1).
Setting $\eta=-v$, we see that \eqref{eq.rho} holds.
Moreover, by Theorem~3.1 in \cite{gk} again, we obtain the
existence of universal $\varepsilon>0$ so that 
$\|\eta\|_2<1$ if $\|\rho\|_2<\varepsilon$.
Thus the proof completes by applying Theorem~\ref{thm.t.q}.
\end{proof}

\section{Linear transformations}
\label{sec.TtoQ}

In this section, we shall see the two-way relationship
between $\mathcal{T}_0$ and $\mathcal{Q}_0$, which has more
explicit representation than that between $\mathcal{T}$ and
$\mathcal{Q}$.

For $\chi,\sigma\in C([0,T];\mathbb{R}^{d\times d})$, define
the linear transformation 
$F_\chi:\mathcal{W}\to \mathcal{W}$ and the Wiener
functional $\mathfrak{p}_\sigma:\mathcal{W}\to \mathbb{R}$
by 
\[
    F_\chi=-\int_0^\bullet \chi(t)\theta(t) dt
    \quad\text{and}\quad
    \mathfrak{p}_\sigma
      =\int_0^T \langle \sigma(t)\theta(t),
        d\theta(t)\rangle.
\]
Defining $\eta_\chi\in \mathcal{S}_2$ by
$\eta_\chi(t,s)=\chi(t)$ for $t>s$ and 
$\eta_\chi(t,t)=(\chi(t)+\chi(t)^\dagger)/2$, we see that
$F_\chi=G_{\eta_\chi}$ and
$\mathfrak{p}_\sigma=\mathfrak{q}_{\eta_\sigma}$. 
Thus $\mathcal{T}_0$ and $\mathcal{Q}_0$ are rewritten as 
\[
    \mathcal{T}_0=\{F_\chi;
       \chi\in C([0,T];\mathbb{R}^{d\times d})\}
    \quad\text{and}\quad
    \mathcal{Q}_0=\{\mathfrak{p}_\sigma;
       \sigma\in C([0,T];\mathbb{R}^{d\times d})\}.
\]

By Theorem~\ref{thm.t.q}, we obtain the
way from $\mathcal{T}_0$ to $\mathcal{Q}_0$.

\begin{theorem}\label{thm.TtoQ}
Suppose that $\chi\in C([0,T];\mathbb{R}^{d\times d})$
satisfies that 
\begin{equation}\label{ass.chi}
    T\|\chi\|_\infty<1, 
\end{equation}
where $\|\chi\|_\infty=\sup_{t\in[0,T]}|\chi(t)|$.
Define $\sigma\in C([0,T];\mathbb{R}^{d\times d})$ by
\begin{equation}\label{eq.sigma}
    \sigma(t)=\chi(t)-\int_t^T 
               \chi(u)^\dagger\chi(u)du,
    \quad t\in[0,T].
\end{equation}
Then $(\iota+F_\chi)^{-1}$ exists and 
is a continuous linear operator from $\mathcal{W}$ to
itself, and it holds that
\begin{equation}\label{t.TtoQ.1}
    \int_{\mathcal{W}} e^{\mathfrak{p}_\sigma} f d\mu
    =e^{\int_0^T [\text{\rm tr}\,(\chi(t)-\sigma(t))]dt/2} 
     \int_{\mathcal{W}} [f\circ(\iota+F_\chi)^{-1}] d\mu,
    \quad f\in C_b(\mathcal{W}).
\end{equation}
\end{theorem}

\begin{proof}
Define $\rho\in \mathcal{S}_2$ by \eqref{eq.rho} with
$\eta=\eta_\chi$.
Then $\rho(t,s)=\sigma(t)$ for $t>s$.
Hence $\mathfrak{p}_\sigma=\mathfrak{q}_\rho$.

By the definitions of $\eta_\chi$ and $\sigma$, it holds
that
\begin{align*}
    & \|\eta_\chi\|_2^2
      =2\int_0^T t|\chi(t)|^2 dt
      =2\int_0^T \Bigl(\int_t^T |\chi(u)|^2 du\Bigr) dt,
    \\
    & \text{\rm tr}\,(\chi(t)-\sigma(t))
      =\int_t^T |\chi(u)|^2 du.
\end{align*}
These identities imply that 
\[
    \|\eta_\chi\|_2\le T\|\chi\|_\infty<1
    \quad\text{and}\quad
    \frac12 \|\eta_\chi\|_2^2
    =\int_0^T [\text{\rm tr}\,(\chi(t)-\sigma(t))]dt.
\]
Notice that the operator norm of the continuous linear
operator $F_\chi:\mathcal{W}\to \mathcal{W}$ is less than
or equal to $T\|\chi\|_\infty<1$.
Hence $(\iota+F_\chi)^{-1}$ exists and is a continuous
linear operator from $\mathcal{W}$ to itself.

The proof completes by applying Theorem~\ref{thm.t.q}
with $\eta_\chi$ and $f\circ (\iota+F_\chi)^{-1}$ for $\eta$ 
and $f$, respectively.
\end{proof}

As was seen in \cite[Lemma~3.1]{t}, if $\chi$ is represented
as $\chi=\alpha^\prime \alpha^{-1}$ for some 
$\alpha\in C^1([0,T];\mathbb{R}^{d\times d})$, then an
explicit expression of $(\iota+F_\chi)^{-1}$ is available. 
If $\|\chi\|_\infty$ is small, such a representation of
$\chi$ is possible as follows.

\begin{proposition}\label{prop.inverse}
Let $\chi\in C([0,T];\mathbb{R}^{d\times d})$.
Suppose that
$T\sqrt{d}\|\chi\|_\infty e^{T\|\chi\|_\infty}<1$.
Define $\alpha\in C^1([0,T];\mathbb{R}^{d\times d})$ to be
the solution to the first order linear ODE
\[
    \alpha^\prime=\chi \alpha,\quad \alpha(T)=I_d.
\]
Then 
\\
{\rm(i)}
$\alpha$ is non-singular, that is, $\det \alpha(t)\ne0$ for
any $t\in[0,T]$,
\\
{\rm(ii)}
the function $\tilde{F}_\chi:\mathcal{W}\to \mathcal{W}$
defined by
\[
    [\tilde{F}_\chi(w)](t)
    =-\alpha(t)\int_0^t (\alpha^{-1})^\prime(s) w(s) ds,
    \quad w\in \mathcal{W},~t\in[0,T],
\]
satisfies that $(\iota+F_\chi)^{-1}=\iota+\tilde{F}_\chi$.
\end{proposition}

\begin{proof}
To see (i), let $\hat{\alpha}=\alpha(T-\cdot)$.
It holds that
\[
    \hat{\alpha}(t)
    =I_d-\int_0^t \chi(T-s)\hat{\alpha}(s)ds,
    \quad t\in[0,T],
\]
which implies that
\[
    |\hat{\alpha}(t)|
    \le \sqrt{d}+\|\chi\|_\infty 
        \int_0^t |\hat{\alpha}(s)|ds,
    \quad t\in[0,T].
\]
By Gronwall's inequality, this yields that
$\|\alpha\|_\infty\le \sqrt{d} e^{T\|\chi\|_\infty }$.
Hence we have that
\[
    |I_d-\alpha(t)|
    =\Bigl|\int_t^T \chi(s)\alpha(s)ds\Bigr|
    \le T\sqrt{d}\|\chi\|_\infty e^{T\|\chi\|_\infty }
    <1.
\]
Thus $\alpha$ is non-singular.

To see (ii), using (i), rewrite $F_\chi$ as
\[
    [F_\chi(w)](t)
    =-\int_0^t \alpha^\prime(s) \alpha^{-1}(s)w(s)ds,
    \quad w\in \mathcal{W},~t\in[0,T].
\]
By the integration by parts on $[0,T]$, 
a direct computation implies that
$(\iota+F_\chi)\circ(\iota+\tilde{F}_\chi)
=(\iota+\tilde{F}_\chi)\circ(\iota+F_\chi)=\iota$.
\end{proof}

Applying Proposition~\ref{prop.inverse}, we have a precise
representation of the conditional expectation
$\mathbb{E}[e^{\mathfrak{p}_\sigma}|\theta(t)=x]$ of
$e^{\mathfrak{p}_\sigma}$ given the condition $\theta(t)=x$.

\begin{proposition}\label{prop.gauss}
Let $\chi$ and $\alpha$ be as in
Proposition~\ref{prop.inverse}.
Define $\sigma$ by \eqref{eq.sigma}.
For $t\in(0,T]$, set
\[
    v_t(\alpha)=\int_0^t (\alpha(t)\alpha(s)^{-1})
         (\alpha(t)\alpha(s)^{-1})^\dagger ds.
\]
Then $v_t(\alpha)$ is positive definite and it holds that
\[
    \mathbb{E}[e^{\mathfrak{p}_\sigma}|\theta(t)=x]
    =e^{\int_0^T[\text{\rm tr}\,(\chi(t)-\sigma(t))]dt/2}
     g_{v_t(\alpha)}(x) \sqrt{2\pi t}^d e^{|x|^2/(2t)},
    \quad x\in \mathbb{R}^d,
\]
where 
\[
    g_{v_t(\alpha)}(x)
    =\frac1{\sqrt{(2\pi)^d \det v_t(\alpha)}}
     e^{-\langle v_t(\alpha)^{-1}x,x\rangle/2},
    \quad x\in \mathbb{R}^d.
\]
\end{proposition}

\begin{proof}
In what follows, we fix $t\in(0,T]$.
It is easy to see that $v_t(\alpha)$ is positive definite. 

By It\^o's formula, we have that
\[
    \int_0^t (\alpha^{-1})^\prime(s)\theta(s)ds
    =\alpha(t)^{-1}\theta(t)
      -\int_0^t \alpha(s)^{-1} d\theta(s).
\]
Hence it holds that
\[
    [\iota+\tilde{F}_\chi](t)
     =\alpha(t)\int_0^t \alpha(s)^{-1} d\theta(s),
\]
where $[\iota+\tilde{F}_\chi](t)$ is the random variable
whose value at $w\in \mathcal{W}$ is 
$[(\iota+\tilde{F}_\chi)(w)](t)$.
Thus it holds that
\[
    \int_{\mathcal{W}} \varphi\bigl(
        [\iota+\tilde{F}_\chi](t)\bigr)
     d\mu
    =\int_{\mathbb{R}^d} \varphi(x) g_{v_t(\alpha)}(x)dx,
    \quad\varphi\in C_b(\mathbb{R}^d).
\]

The assumption that
$T\sqrt{d}\|\chi\|_\infty e^{T\|\chi\|_\infty}<1$ yields that
$T\|\chi\|_\infty<1$.
By Theorem~\ref{thm.TtoQ} and
Proposition~\ref{prop.inverse}, 
the identity \eqref{t.TtoQ.1} holds with
$\iota+\tilde{F}_\chi$ for $(\iota+F_\chi)^{-1}$.
Hence we have that
\begin{align*}
    & \int_{\mathbb{R}^d}
        \mathbb{E}[e^{\mathfrak{p}_\sigma}|\theta(t)=x] 
        \varphi(x) \frac1{\sqrt{2\pi t}^d}
        e^{-|x|^2/(2t)} dx
     = \int_{\mathcal{W}} e^{\mathfrak{p}_\sigma}
           \varphi(\theta(t))d\mu
    \\
    & \quad
     = e^{\int_0^T[\text{\rm tr}\,(\chi(t)-\sigma(t))]dt/2}
        \int_{\mathbb{R}^d} \varphi(x) g_{v_t(\alpha)}(x)dx,
    \qquad \varphi\in C_b(\mathbb{R}^d).
\end{align*}
This completes the proof.
\end{proof}

We now proceed to showing the way from $\mathcal{Q}_0$ to 
$\mathcal{T}_0$.
Introduce the conditions on $\varepsilon,\delta>0$ such that 
\begin{align}
    & 2\varepsilon T\sqrt{d} 
      \{1+T\sqrt{d}(1+\varepsilon)\}
      e^{T(\sqrt{d}+2 \varepsilon+\varepsilon^2)}<1,
 \label{eq.eps.1}
    \\
    & \varepsilon T\{
      1+T\sqrt{d} K_0 (1+\varepsilon)
      e^{T(\sqrt{d}+2 \varepsilon+\varepsilon^2)}\}
      <1,
 \label{eq.eps.2}
    \\
    & \delta T(2\sqrt{d}\vee K_0)
          \{1+T(\sqrt{d}+\delta)
          e^{T(\sqrt{d}+\delta)}\}<1,
 \label{eq.delta}
\end{align}
where 
\[
    K_0=\sup\bigl\{|M^{-1}|;M\in \mathbb{R}^{d\times d},
        |M-I_d|<\tfrac12\bigr\}.    
\]
Put $\boldsymbol{\epsilon}(\sigma)=\|\sigma\|_\infty$ for  
$\sigma\in C([0,T];\mathbb{R}^{d\times d})$
and 
$\boldsymbol{\delta}(\sigma)
  =|\sigma(T)|+\|\sigma^\prime\|_\infty
    +2\|\sigma_A\|_\infty$
for $\sigma\in C^1([0,T];\mathbb{R}^{d\times d})$.
Our second goal of this section is the following.

\begin{theorem}\label{thm.QtoT}
Suppose that $\varepsilon>0$ and $\delta>0$ satisfy
\eqref{eq.eps.1}, \eqref{eq.eps.2}, and \eqref{eq.delta}.
\\
{\rm(i)}
Let $\sigma\in C([0,T];\mathbb{R}^{d\times d})$.
Suppose that 
$\boldsymbol{\epsilon}(\sigma)<\varepsilon$.
Then the following assertions hold.

{\rm(a)}
There exists $S\in C^1([0,T];\mathbb{R}^{d\times d})$
obeying the ODE
\begin{equation}\label{eq.riccati}
   S^\prime=-S^2-\sigma^\dagger S-S\sigma
        -\sigma^\dagger \sigma,
   \quad  S(T)=0.
\end{equation}

{\rm(b)}
The function $\chi=S+\sigma$ satisfies \eqref{ass.chi} and 
\eqref{eq.sigma}, and it holds that
\begin{equation}\label{t.q.l.A.2}
    \int_{\mathcal{W}} e^{\mathfrak{p}_\sigma} f d\mu
    =e^{\int_0^T [\text{\rm tr}\, S(t)] dt/2}
      \int_{\mathcal{W}} [f\circ(\iota+F_\chi)^{-1}] d\mu,
    \quad f\in C_b(\mathcal{W}).
\end{equation}

\smallskip
\noindent
{\rm(ii)}
Let $\sigma\in C^1([0,T];\mathbb{R}^{d\times d})$.
Suppose that $\boldsymbol{\delta}(\sigma)<\delta$.
Then the following assertions hold.

{\rm(a)}
The solution $A\in C^2([0,T];\mathbb{R}^{d\times d})$
to the ODE
\begin{equation}\label{eq.jacobi}
    A^{\prime\prime}-2\sigma_A A^\prime
      -\sigma^\prime A=0, \quad
    A(T)=I_d,~A^\prime(T)=\sigma(T)
\end{equation}
is non-singular, that is, 
$\det A(t)\ne0$ for any $t\in[0,T]$.

{\rm(b)}
The function $\chi=A^\prime A^{-1}$ satisfies
\eqref{ass.chi} and \eqref{eq.sigma}, and it holds that
\begin{equation}\label{t.q.l.B.2}
    \int_{\mathcal{W}} e^{\mathfrak{p}_\sigma} f d\mu
    =\frac{e^{-\int_0^T [\text{\rm tr}\,\sigma_S(t)] dt/2}}{
           \sqrt{\det A(0)}}
      \int_{\mathcal{W}} [f\circ(\iota+F_\chi)^{-1}] d\mu,
    \quad f\in C_b(\mathcal{W}),
\end{equation}
where $\sigma_S=\frac12(\sigma+\sigma^\dagger)$.
\end{theorem}

The proof of the theorem is broken into several steps, each
step being a lemma. 
We start with an elementary lemma on linear ODEs.

\begin{lemma}\label{lem.ode}
Let $\xi_1,\xi_2\in \mathbb{R}^{d\times d}$ and
$\gamma_{ij}\in C([0,T];\mathbb{R}^{d\times d})$,
$i,j=1,2$.
Define 
$\phi_1,\phi_2\in C^1([0,T];\mathbb{R}^{d\times d})$ by 
the ODE on $\mathbb{R}^{2d\times d}$;
\begin{equation}\label{l.ode.1}
    \begin{pmatrix} \phi_1 \\ \phi_2\end{pmatrix}^\prime
    =\begin{pmatrix} 
          \gamma_{11} & \gamma_{12} \\ 
          \gamma_{21} & \gamma_{22} 
      \end{pmatrix}
      \begin{pmatrix} \phi_1 \\ \phi_2\end{pmatrix},
    \quad
    \begin{pmatrix} \phi_1(0) \\ \phi_2(0)\end{pmatrix}
    =\begin{pmatrix} \xi_1 \\ \xi_2\end{pmatrix}.
\end{equation}
Then it holds that
\begin{align*}
    \|\phi_2\|_\infty
      \le & |\xi_2|
           +T \bigl({\textstyle\sum_j}|\xi_j|\bigr)
            \bigl({\textstyle\sum_j}
            \|\gamma_{2j}\|_\infty\bigr)
          e^{T\sum_{i,j}\|\gamma_{ij}\|_\infty},
    \\
    \|\phi_1-\xi_1\|_\infty
      \le & T\|\gamma_{11}\|_\infty 
          \bigl({\textstyle\sum_j}|\xi_j|\bigr)
          e^{T\sum_{i,j}\|\gamma_{ij}\|_\infty}
    \\
    & +T\|\gamma_{12}\|_\infty
          \Bigl\{|\xi_2|+T
           \bigl({\textstyle\sum_j}|\xi_j|\bigr)
           \bigl({\textstyle\sum_j}
             \|\gamma_{2j}\|_\infty\bigr)
          e^{T\sum_{i,j}\|\gamma_{ij}\|_\infty} \Bigr\},
\end{align*}
where $\sum_j$ and $\sum_{i,j}$ are the abbreviations of
$\sum_{j=1}^2$ and $\sum_{i,j=1}^2$, respectively.
Moreover, if $\det\phi_1(t)\ne0$ for any $t\in[0,T]$, then
$\psi=\phi_2\phi_1^{-1}$ obeys the ODE
\[
    \psi^\prime
    =-\psi\gamma_{12}\psi+\gamma_{22}\psi
     -\psi\gamma_{11}+\gamma_{21},\quad
    \psi(0)=\xi_2\xi_1^{-1}.
\]
\end{lemma}

\begin{proof}
The last assertion is easily shown, and it is a well-known 
method to solve matrix Riccati ODEs (cf.\cite{f}). 

Taking the sum of the norms of upper and lower halves of 
\eqref{l.ode.1}, we have that
\[
    {\textstyle\sum_j}|\phi_j(t)|
    \le {\textstyle\sum_j}|\xi_j|
       +\bigl({\textstyle\sum_{i,j}}
          \|\gamma_{ij}\|_\infty \bigr)
        \int_0^t \bigl({\textstyle\sum_j}|\phi_j(s)| 
              \bigr) ds,
    \quad t\in[0,T].
\]
Applying Gronwall's inequality, we obtain that
\begin{equation}\label{l.ode.21}
    \|\phi_j\|_\infty
    \le \bigl({\textstyle\sum_j}|\xi_j|\bigr)
          e^{T\sum_{i,j}\|\gamma_{ij}\|_\infty},
    \quad j=1,2.
\end{equation}
Substitute this into the lower half of \eqref{l.ode.1}, we
obtain the first inequality.
Plugging the first inequality and \eqref{l.ode.21} for $j=1$
into the upper half of \eqref{l.ode.1}, we arrive at the
second inequality.
\end{proof}

We now proceed to the proof of the assertion (i) of
Theorem~\ref{thm.QtoT}.
In the following two lemmas, we always assume that
$\sigma\in C([0,T];\mathbb{R}^{d\times d})$ and 
$\boldsymbol{\epsilon}(\sigma)<\varepsilon$.

\begin{lemma}\label{lem.riccati}
There exists $S\in C^1([0,T];\mathbb{R}^{d\times d})$
obeying the ODE \eqref{eq.riccati}.
\end{lemma}

\begin{proof}
For $\kappa\in C([0,T];\mathbb{R}^{d\times d})$, define
$\hat{\kappa}\in C([0,T];\mathbb{R}^{d\times d})$ by
$\hat{\kappa}(t)=\kappa(T-t)$, $t\in[0,T]$.
Then the ODE \eqref{eq.riccati} to be solved turns into
\begin{equation}\label{t.q.l.A.22}
    \hat{S}^\prime
    =\hat{S}^2+\hat{\sigma}^\dagger \hat{S}
     +\hat{S} \hat{\sigma}
     +\hat{\sigma}^\dagger \hat{\sigma},
    \quad
    \hat{S}(0)=0.
\end{equation}
Define 
$\begin{pmatrix} \phi_1 \\ \phi_2\end{pmatrix}
  \in C^1([0,T];\mathbb{R}^{2d\times d})$ 
by the ODE
\[
    \begin{pmatrix} \phi_1 \\ \phi_2\end{pmatrix}^\prime
    =\begin{pmatrix}
       -\hat{\sigma} & -I_d \\
       \hat{\sigma}^\dagger\hat{\sigma} & \hat{\sigma}^\dagger
     \end{pmatrix}
     \begin{pmatrix} \phi_1 \\ \phi_2\end{pmatrix},
    \quad 
     \begin{pmatrix} \phi_1(0) \\ \phi_2(0)\end{pmatrix} 
     =\begin{pmatrix} I_d \\ 0\end{pmatrix}.    
\]
The second inequality in Lemma~\ref{lem.ode} yields that
\[
    \|\phi_1-I_d\|_\infty
    \le \boldsymbol{\epsilon}(\sigma) T\sqrt{d} 
       \{1+T\sqrt{d}(1+\boldsymbol{\epsilon}(\sigma))\}
      e^{T(\sqrt{d}+2 \boldsymbol{\epsilon}(\sigma)
                +\boldsymbol{\epsilon}(\sigma)^2)}.
\]
Hence, by \eqref{eq.eps.1},
$\|\phi_1-I_d\|_\infty<\frac12$,
and $\det\phi_1(t)\ne0$ for any $t\in[0,T]$.
Due to the same lemma, we see that the function
$\phi_2\phi_1^{-1}$ solves the ODE \eqref{t.q.l.A.22}.
Thus 
$S\in C^1([0,T];\mathbb{R}^{d\times d})$, determined by the
relation that $\hat{S}=\phi_2\phi_1^{-1}$, is the solution
to the ODE \eqref{eq.riccati}. 
\end{proof}

\begin{lemma}\label{lem.S.chi}
Let $S\in C^1([0,T];\mathbb{R}^{d\times d})$ be as in
Lemma~\ref{lem.riccati}.
Then $\chi=S+\sigma$ satisfies \eqref{ass.chi} and 
\eqref{eq.sigma}, and \eqref{t.q.l.A.2} holds.
In particular, the assertion {\rm(i)} of
Theorem~\ref{thm.QtoT} holds.
\end{lemma}

\begin{proof}
We first show that $\chi=S+\sigma$ satisfies
\eqref{eq.sigma}. 
Since $S$ obeys the Riccati ODE \eqref{eq.riccati}, it holds 
that
\begin{equation}\label{t.q.l.A.21}
    S(t)-\int_t^T (S(s)+\sigma(s)^\dagger)
        (S(s)+\sigma(s))ds=0.
\end{equation}
Taking the transpose of this identity, we see that
$S^\dagger=S$.
Hence $\chi^\dagger=S+\sigma^\dagger$.
Substituting this into \eqref{t.q.l.A.21},
and adding $\sigma$ to both sides of the resulting
identity, we see that the identity \eqref{eq.sigma} holds.

We next show that $\chi=S+\sigma$ satisfies
\eqref{ass.chi}. 
Since $\|\phi_1-I_d\|_\infty<\frac12$ as was seen in the
proof of the previous lemma, the first inequality in
Lemma~\ref{lem.ode} implies that 
\[
    \|S\|_\infty
    =\|\phi_2\phi_1^{-1}\|_\infty
    \le \boldsymbol{\epsilon}(\sigma) T
        \sqrt{d} K_0 (1+\boldsymbol{\epsilon}(\sigma))
        e^{T(\sqrt{d}+2 \boldsymbol{\epsilon}(\sigma)
          +\boldsymbol{\epsilon}(\sigma)^2)},
\]
where $\phi_1,\phi_2$ are the functions given in the proof
of Lemma~\ref{lem.riccati} to construct $S$.
By \eqref{eq.eps.2}, this implies that
$T\|S+\sigma\|_\infty<1$, and hence 
$\chi=S+\sigma$ satisfies \eqref{ass.chi}.

Since $\chi-\sigma=S$, the identity \eqref{t.q.l.A.2}
follows from Theorem~\ref{thm.TtoQ}.
\end{proof}

We now give the proof of the assertion (ii) of
Theorem~\ref{thm.QtoT}.
In the following two lemmas, we always assume that
$\sigma\in C^1([0,T];\mathbb{R}^{d\times d})$ and
$\boldsymbol{\delta}(\sigma)<\delta$.

\begin{lemma}\label{lem.jacobi}
The solution $A\in C^2([0,T];\mathbb{R}^{d\times d})$ to the
ODE \eqref{eq.jacobi} is non-singular.
\end{lemma}

\begin{proof}
Define $\phi_1,\phi_2\in C([0,T];\mathbb{R}^{d\times d})$ by
$\phi_1=A(T-\cdot)$ and $\phi_2=-A^\prime(T-\cdot)$.
It then holds that
\[
    \begin{pmatrix} \phi_1 \\ \phi_2\end{pmatrix}^\prime
    =\begin{pmatrix}
       0 & I_d \\
       \sigma^\prime(T-\cdot) & -2\sigma_A(T-\cdot)
     \end{pmatrix}
     \begin{pmatrix} \phi_1 \\ \phi_2\end{pmatrix},
    \quad 
     \begin{pmatrix} \phi_1(0) \\ \phi_2(0)\end{pmatrix} 
     =\begin{pmatrix} I_d \\ -\sigma(T)\end{pmatrix}.    
\]
The second inequality in Lemma~\ref{lem.ode} yields that
\[
    \|A-I_d\|_\infty
    =\|\phi_1-I_d\|_\infty
    \le \boldsymbol{\delta}(\sigma) T\sqrt{d}
        \{1+T(\sqrt{d}+\boldsymbol{\delta}(\sigma))
             e^{T(\sqrt{d}+\boldsymbol{\delta}(\sigma))}\}.
\]
By \eqref{eq.delta},
$\|A-I_d\|_\infty<\frac12$, and hence $A$ is non-singular.
\end{proof}

\begin{lemma}\label{lem.A.jacobi}
Let $A$ be as in Lemma~\ref{lem.jacobi}.
Then $\chi=A^\prime A^{-1}$ satisfies \eqref{ass.chi} and 
\eqref{eq.sigma}, and \eqref{t.q.l.B.2} holds.
In particular, the assertion {\rm (ii)} of
Theorem~\ref{thm.QtoT}.
holds.
\end{lemma}

\begin{proof}
We first show that $\chi=A^\prime A^{-1}$ satisfies 
\eqref{eq.sigma}.
To do so, put $S=\chi-\sigma$.
Then $\chi$ and $S$ are both in 
$C^1([0,T];\mathbb{R}^{d\times d})$ and obey the following ODEs:
\begin{align}
    & \chi^\prime=-\chi^2+2\sigma_A\chi+\sigma^\prime,
      && \chi(T)=\sigma(T),
 \label{t.q.l.B.22}
    \\
    & S^\prime=-S^2-\sigma^\dagger S-S\sigma
           -\sigma^\dagger\sigma,
      && S(T)=0.
 \nonumber
\end{align}
Since $S^\dagger$ solves the same ODE as $S$ does,
$S^\dagger=S$.
Hence $\chi^\dagger=\chi-2\sigma_A$.
Plugging this into \eqref{t.q.l.B.22}, we obtain that
$\chi^\prime+\chi^\dagger\chi=\sigma^\prime$ and 
$\chi(T)=\sigma(T)$.
Thus \eqref{eq.sigma} holds.

We next see that $\chi=A^\prime A^{-1}$ satisfies
\eqref{ass.chi}. 
Due to the first inequality in Lemma~\ref{lem.ode}
and \eqref{eq.delta}, we have that
\[
    \|A^\prime\|_\infty
    =\|\phi_2\|_\infty
    \le \boldsymbol{\delta}(\sigma)
          \{1+T(\sqrt{d}+\boldsymbol{\delta}(\sigma))
        e^{T(\sqrt{d}+\boldsymbol{\delta}(\sigma))}\}
    <\frac{1}{TK_0}.
\]
As was seen in the previous proof, it holds that
$\|A-I_d\|_\infty<\frac12$, and hence 
\[
    T\|\chi\|_\infty \le TK_0\|A^\prime\|_\infty<1.
\]
Thus \eqref{ass.chi} holds.

We finally show the identity \eqref{t.q.l.B.2}.
Remember that the mapping $t\mapsto\det A(t)$ obeys the 
ODE
\[
    (\det A)^\prime =[\text{\rm tr}\,(A^\prime A^{-1})]
       \det A,\quad
    \det A(T) =1.
\]
Due to the definition of $\chi$, this implies that
\[
    \det A(t)
    =e^{-\int_t^T \text{\rm tr}\,[A^\prime(s)A^{-1}(s)]ds}
    =e^{-\int_t^T \text{\rm tr}\,\chi(s) ds}.
\]
Since $\text{\rm tr}\,\sigma=\text{\rm tr}\,\sigma_S$, in
conjunction with Theorem~\ref{thm.TtoQ}, this yields
\eqref{t.q.l.B.2}. 
\end{proof}

\begin{remark}\label{rem.3.1}
Since 
$\boldsymbol{\epsilon}(\sigma)\le 
 |\sigma(T)|+T\|\sigma^\prime\|_\infty$,
it holds that
$\boldsymbol{\epsilon}(\sigma)\le 
 (1+T)\boldsymbol{\delta}(\sigma)$.
Thus, if $\delta>0$ in \eqref{eq.delta} is chosen so that 
$(1+T)\delta<\varepsilon$, then the assertions (i) and (ii)
of Theorem~\ref{thm.QtoT} are both applicable. 
Further, in this case, the Riccati ODE \eqref{eq.riccati} 
follows from the ODE \eqref{eq.jacobi}.
In fact, let $A\in C^2([0,T];\mathbb{R}^{d\times d})$  be
the solution to the linear ODE \eqref{eq.jacobi}, and set
$\chi=A^\prime A^{-1}$ and $S=\chi-\sigma$.
As was seen just after \eqref{t.q.l.B.22}, $S$ obeys the ODE
\eqref{eq.riccati}. 
\end{remark}

Since 
$\boldsymbol{\epsilon}(\lambda\sigma)
   =|\lambda|\boldsymbol{\epsilon}(\sigma)$ and
$\boldsymbol{\delta}(\lambda\sigma)
 =|\lambda|\boldsymbol{\delta}(\sigma)$ for
$\lambda\in \mathbb{R}$, the previous theorem implies the
following.

\begin{corollary}\label{cor.q.l}
Suppose that $\varepsilon>0$ and $\delta>0$ satisfy
\eqref{eq.eps.1}, \eqref{eq.eps.2}, and \eqref{eq.delta}.
\\
{\rm(i)}
Let $\sigma\in C([0,T];\mathbb{R}^{d\times d})$.
Suppose that $\lambda\in \mathbb{R}$ satisfies that
$|\lambda|\boldsymbol{\epsilon}(\sigma)<\varepsilon$.
Then the following assertions hold.

{\rm(a)}
There exists 
$S_\lambda\in C^1([0,T];\mathbb{R}^{d\times d})$ obeying the
ODE 
\[
   S_\lambda^\prime
      =-S_\lambda^2 -\lambda \sigma^\dagger S_\lambda
       -\lambda S_\lambda\sigma
      -\lambda^2 \sigma^\dagger\sigma,
   \quad S_\lambda(T)=0.
\]

{\rm(b)}
Let $\chi_\lambda=S_\lambda+\lambda\sigma$.
Then it holds that
\[
    \int_{\mathcal{W}} e^{\lambda \mathfrak{p}_\sigma} f d\mu
    =e^{\int_0^T [\text{\rm tr}\, S_\lambda(t)] dt/2}\int_{\mathcal{W}} 
      [f\circ(\iota+F_{\chi_\lambda})^{-1}] d\mu,
    \quad f\in C_b(\mathcal{W}).
\]

\smallskip\noindent
{\rm(ii)}
Let $\sigma\in C^1([0,T];\mathbb{R}^{d\times d})$.
Suppose that 
$\lambda\in \mathbb{R}$ satisfies that
$|\lambda|\boldsymbol{\delta}(\sigma)<\delta$.
Then the following assertions hold.

{\rm(a)}
The solution 
$A_\lambda \in C^2([0,T];\mathbb{R}^{d\times d})$ to 
the ODE
\[
    A_\lambda^{\prime\prime}
     -2\lambda\sigma_A A_\lambda^\prime
     -\lambda\sigma^\prime A_\lambda=0, \quad
    A_\lambda(T)=I_d,~A_\lambda^\prime(T)=\lambda\sigma(T)
\]
is non-singular, that is,
$\det A_\lambda(t)\ne0$ for any $t\in[0,T]$.

{\rm(b)}
Let $\chi_\lambda=A_\lambda^\prime A_\lambda^{-1}$.
Then it holds that
\[
    \int_{\mathcal{W}} e^{\lambda \mathfrak{p}_\sigma} f d\mu
    =\frac{e^{-\lambda \int_0^T [\text{\rm tr}\,\sigma_S(t)] dt/2}}{
           \sqrt{\det A_\lambda(0)}}
      \int_{\mathcal{W}} 
      [f\circ(\iota+F_{\chi_\lambda})^{-1}] d\mu,
     ~~ f\in C_b(\mathcal{W}).
\]
\end{corollary}

In the remaining of this section, we consider the Wiener
functional $\mathfrak{q}:\mathcal{W}\to\mathbb{R}$ given by 
\[
    \mathfrak{q}
    =\int_0^T \langle\gamma(t)\theta(t),d\theta(t)
      \rangle
    +\frac12 \int_0^T \langle\kappa(t)\theta(t),\theta(t)
       \rangle dt,
\]
where $\gamma,\kappa\in C([0,T];\mathbb{R}^{d\times d})$.
By It\^o's formula, we have that
\begin{equation}\label{eq.q}
    \mathfrak{q}=\mathfrak{p}_\sigma
         +\frac12 \int_0^T\Bigl(
          \int_t^T[\text{\rm tr}\,\kappa_S(s)]ds \Bigr)dt,
\end{equation}
where $\kappa_S=\frac12(\kappa+\kappa^\dagger)$, and 
$\sigma\in C([0,T];\mathbb{R}^{d\times d})$ is given by
\[
    \sigma(t)=\gamma(t)+\int_t^T \kappa_S(s)ds,
    \quad t\in[0,T].
\]
As an application of Theorem~\ref{thm.QtoT}, we have the
following.

\begin{corollary}\label{cor.quad}
Suppose that $\varepsilon,\delta>0$ satisfy
\eqref{eq.eps.1}, \eqref{eq.eps.2}, and \eqref{eq.delta}.
Let $\gamma,\kappa,\sigma,\mathfrak{q}$ be as above.

\noindent
{\rm(i)}
Suppose that 
$\|\gamma\|_\infty+T\|\kappa\|_\infty<\varepsilon$.
Then the following assertions hold.

{\rm(a)}
There exists $S\in C^1([0,T];\mathbb{R}^{d\times d})$
obeying the ODE
\[
   S^\prime=-S^2-\sigma^\dagger S-S\sigma
        -\sigma^\dagger \sigma,
   \quad  S(T)=0.
\]

{\rm(b)}
Let $\chi=S+\sigma$. 
Then it holds that
\[
    \int_{\mathcal{W}} e^{\mathfrak{q}} f d\mu
    =e^{\int_0^T [\text{\rm tr}\,(S(t)+\int_t^T\kappa_S(s)ds)] dt/2}
      \int_{\mathcal{W}} [f\circ(\iota+F_\chi)^{-1}] d\mu,
    \quad f\in C_b(\mathcal{W}).
\]

\smallskip\noindent
{\rm(ii)}
Suppose that 
$\gamma\in C^1([0,T];\mathbb{R}^{d\times d})$ and 
$|\gamma(T)|+\|\gamma^\prime-\kappa\|_\infty
  +2\|\gamma_A\|_\infty<\delta$, where
$\gamma_A=\frac12(\gamma-\gamma^\dagger)$.
Then the following assertions hold.

{\rm(a)}
The solution $A\in C^2([0,T];\mathbb{R}^{d\times d})$
to the ODE
\[
    A^{\prime\prime}-2\gamma_A A^\prime
      +(\kappa_S-\gamma^\prime) A=0, \quad
    A(T)=I_d,~A^\prime(T)=\gamma(T)
\]
is non-singular.

{\rm(b)}
Let $\chi=A^\prime A^{-1}$.
Then it holds that
\[
    \int_{\mathcal{W}} e^{\mathfrak{q}} f d\mu
    =\frac{e^{-\int_0^T [\text{\rm tr}\,\gamma_S(t)] dt/2}}{
           \sqrt{\det A(0)}}
      \int_{\mathcal{W}} [f\circ(\iota+F_\chi)^{-1}] d\mu,
    \quad f\in C_b(\mathcal{W}),
\]
where
$\gamma_S=\frac12(\gamma+\gamma^\dagger)$.
\end{corollary}

\begin{proof}
Notice that
$\|\sigma\|_\infty\le \|\gamma\|_\infty+T\|\kappa\|_\infty$,
$\sigma(T)=\gamma(T)$, 
$\sigma^\prime=\gamma^\prime-\kappa_S$,
$\sigma_S=\gamma_S+\int_\bullet^T \kappa_S(s)ds$,
$\sigma_A=\gamma_A$, and
\[
    \int_{\mathcal{W}} e^{\mathfrak{q}} f d\mu
    =e^{\int_0^T[\text{\rm tr}\,(\int_t^T\kappa_S(s)ds)]dt/2} 
     \int_{\mathcal{W}} e^{\mathfrak{p}_\sigma} f d\mu.
\]
The assertions follows by applying Theorem~\ref{thm.QtoT} to
$\mathfrak{p}_\sigma$.
\end{proof}

The evaluation of 
$\int_{\mathcal{W}} e^{\mathfrak{q}} f d\mu$ as stated in
Corollary~\ref{cor.quad}(ii)(b) was first pointed out by
Cameron and Martin \cite{cm,cm1} when $d=1$.
In their case, $\gamma=0$ and $\mathfrak{q}$ is the weighted
square of sample norm  
$\int_0^T \kappa(t)|\theta(t)|^2 dt$.
The corresponding ODE is the Sturm-Liouville equation 
\[
    f^{\prime\prime}+\kappa f=0,\quad f(T)=1,~f^\prime(T)=0.
\]
If $\kappa\equiv1$, then it corresponds to the harmonic
oscillator (\cite{cm1,kac}, also see
\cite[Subsection~5.8.1]{mt}). 
When $d=2$, 
$\gamma=\begin{pmatrix} 0 & -1 \\ 1 & 0 \end{pmatrix}$, and 
$\kappa=0$, $\mathfrak{q}$ is L\'evy's stochastic area and
the evaluation presents L\'evy's stochastic area formula
(\cite{levy,yor}, also see \cite[Subsection~5.8.2]{mt}). 
Such an evaluation was extended to general dimensions 
by the author \cite{t} with the additional assumption that 
$\gamma^\dagger=-\gamma$.
The extension was made by using the Girsanov formula and it
was applied in \cite{it} to representing heat kernels of
step-two nilpotent Lie groups.


\bigskip
\hfill
\begin{minipage}{285pt}
Setsuo Taniguchi
\\
{\tt
E-mail:taniguchi.setsuo.772@m.kyushu-u.ac.jp}
\end{minipage}

\end{document}